\documentclass[12pt]{article}
\usepackage{amssymb}
\usepackage{amsmath}
\usepackage{cases}
\usepackage{version}
\usepackage{graphicx}

\newtheorem{theor}{Theorem}[section]

\newcommand{\cqfd}{\hfill $\square$}
\newcommand{\R}{\mathbb R}

\def \E{\mathbb E}
\newcommand{\ind}{independent }
\newcommand{\crt}{characterization }

\newcommand{\egd}{=_d}
\title{Characterizations of GIG laws: a survey complemented with two new results
}
\author{Angelo Efoevi Koudou \ and \ Christophe Ley \\ {\it \small Universit\'e de Lorraine and CNRS, Institut Elie Cartan de Lorraine}  \\
 {\it \small ECARES and D\' epartement de Math\' ematique, Universit\'{e} libre de Bruxelles}} 
\date{}
\begin{document}
\maketitle
\markboth{ \small GIG characterizations}{ \small GIG characterizations}

\begin{abstract}
Several characterizations of the Generalized Inverse Gaussian (GIG) distribution have been proposed in the literature, especially over the past two decades. These characterization theorems are  surveyed, and two new characterizations are established, one based on maximum likelihood estimation  and the other is a Stein characterization.
\end{abstract}


\section{Introduction}

The \emph{Generalized Inverse Gaussian} (hereafter GIG) distribution has been proposed by Good~\cite{goo} in his study of population frequencies, yet its first appearance can be traced back to Etienne Halphen in the forties~\cite{Halp}, whence the GIG is sometimes called Halphen Type A distribution. This distribution with parameters $p\in\R$, $a>0$, $b>0$  has density 
\begin{equation*}
f_{p,a,b}(x):=\frac{(a/b)^{p/2}}{2K_p(\sqrt{ab})} x^{p-1}e^{-(ax+b/x)/2}, \quad x>0,
\end{equation*}
with $K_p$ the modified Bessel function of the third kind. The parameters $a$ and $b$ regulate both the concentration and scaling of the densities, the former via $\sqrt{ab}$ and the latter via $\sqrt{b/a}$. This is the reason why some authors (e.g., J\o rgensen~\cite{jor}) introduce the parameters $\theta=\sqrt{ab}$ and $\eta=\sqrt{b/a}$, leading to a GIG density of the form  
\begin{equation}\label{GIG2}
f_{p,\theta,\eta}(x):= \frac{1}{ 2\eta K_p (\theta)} \left(\frac{x}{\eta}\right)^{p -1}
e^{-\theta( x/\eta+\eta /x)/2},  \quad x>0.
\end{equation}
The parameter $p$ bears no concrete statistical meaning, but some particular values of $p$ lead to well-known sub-models of the GIG such as the Inverse Gaussian  ($p=-\frac{1}{2}$), the  {Reciprocal Inverse Gaussian} ($p=\frac{1}{2}$), the  {hyperbolic} ($p=0$, in which case one also speaks of the Harmonic law), the positive hyperbolic ($p=1$), the Gamma ($p>0$ and $b=0$) or the Reciprocal Gamma distribution ($p<0$ and $a=0$).

Several papers have investigated the probabilistic  properties of the GIG distributions. For instance, Barndorff-Nielsen and Halgreen~\cite{bar} have studied their convolution properties and  established their infinite divisibility as well as interesting properties such as the equivalence $X \sim {\rm GIG}(p , a, b)\Longleftrightarrow \frac{1}{X} \sim {\rm GIG}(-p , b, a)$. Vallois~\cite{val} has shown that, depending on the sign of $p$, the GIG laws  can be viewed as the distributions of either first or last exit times for certain diffusion processes, generalizing the well-known fact that  the Inverse Gaussian and Reciprocal Inverse Gaussian distributions are respectively the distributions of the first and the last hitting time for a Brownian motion. Madan, Roynette and Yor~\cite{mad} have shown that  the Black-Scholes formula in finance can be expressed in terms of the distribution function of GIG variables. GIG distributions belong to the family of  generalized $\Gamma$-convolutions, as shown in  Halgreen~\cite{hal}  (see also Eberlein and Hammerstein~\cite{ebe} for a detailed proof).  Sichel (\cite{sic1}, ~\cite{sic2}) used this distribution to construct mixtures of Poisson distributions.
The interesting properties of the GIG have led to the definition of matrix GIG distributions, containing the special case of Wishart matrices for $b=0$. For an overview of results on matrix GIG laws, we refer the reader to  Massam and Weso\l owski (\cite{masw1}, ~\cite{masw2}).

The GIG distributions have been used in the modelization of diverse real phenomena such as, for instance, waiting time (J\o rgensen~\cite{jor}), neural activity (Iyengar and Liao~\cite{Iyengar}),  or, most importantly, hydrologic extreme events (see Chebana \emph{et al.}~\cite{Chebana} and references therein). Despite this popularity of GIG models, statistical aspects of the GIG distributions have received much less attention in the literature than their probabilistic counterparts. The cornerstone reference in this respect is J\o rgensen~\cite{jor}, complemented by the stream of literature on general Halphen distributions (of which the GIG or Halphen A is one of three representatives), see e.g. Perreault \emph{et al.}~(\cite{Pera}, \cite{Perb})  or Chebana \emph{et al.}~\cite{Chebana}. Quite recently, Koudou and Ley~\cite{KL} have applied the Le Cam methodology (see~\cite{LeCam}) to  GIG distributions in order to construct optimal (in the maximin sense) testing procedures within the GIG family of distributions.

The present paper is concerned with characterization results of GIG laws. In probability and statistics, a characterization theorem occurs when a given distribution  is the only one which satisfies a certain property. Besides their evident mathematical interest \emph{per se}, characterization theorems also deepen our understanding of the distributions under investigation and sometimes open unexpected paths to innovations which might have been uncovered otherwise. Over the years, an important number of such characterization theorems for GIG laws have been presented in the literature, especially over the past two decades. We therefore propose here a survey of existing characterizations (without pretending to be exhaustive), amended by two novel characterizations, one based on the maximum likelihood estimator (MLE) of the scale parameter $\eta=\sqrt{b/a}$ of GIG laws  and the other is a so-called Stein characterization. 

%
%
%
The paper is organized as follows. A review of known characterizations of GIG distributions is presented in Section 2, while the short Section 3 is devoted to our two new characterizations.

\section{Known characterizations of GIG laws}
\subsection{A characterization by continued fractions}
The following theorem characterizes the GIG distribution as the law of a continued fraction with independent Gamma entries. We adopt the notation $\Gamma(p,a/2)={\rm GIG}(p,a,0)$  for the Gamma distribution with parameters $p,a>0$ and $\egd$ for equality in distribution.
\begin{theor}(Letac and Seshadri,~\cite{lets}).
\begin{itemize}
\item Let $X$ and $Y$ be two \ind random variables  such that $X>0$ and $Y\sim \Gamma(p, a/2)$ for $p, a>0$. Then $X\egd \frac{1}{Y+X}$ if and only if $X \sim {\rm GIG} (-p, a, a)$.
\item Let $X$, $Y_1$ and $Y_2$ be three \ind random variables  such that $X>0$, $Y_1\sim \Gamma(p, b/2)$ and $Y_2\sim \Gamma(p, a/2)$ for $p, a,b>0$. Then $X\egd \frac{1}{Y_1+\frac{1}{Y_2+X}}$ if and only if $X \sim {\rm GIG} (-p, a, b)$.

\item If $(Y_i)_{i\geq
1}$ is a sequence of independent random variables such that
$$\mathcal{L}(Y_{2i-1})=\mathcal{L}(Y_1)=\Gamma(p,b/2)\ \ and\ \
\mathcal{L}(Y_{2i})=\mathcal{L}(Y_2)=\Gamma(p,a/2);\ \ i\geq 1,$$ then
\begin{equation} \label{contf}
\mathcal{L}\left(\displaystyle\frac{1}{Y_1+\displaystyle\frac{1}{Y_2+\displaystyle\frac{1}{Y_3+\ddots}}}\right)=GIG(-p,a,b).\end{equation}
\end{itemize}
\end{theor}\vspace{2mm}

The proof relies on properties of continued fractions, which are exploited to show that the GIG distribution is the unique stationary distribution of the Markov chain $(X_m)_m$ defined by 
$$X_{m+1}=\frac{1}{Y_{m+1}+X_m},$$ in the case $a=b$, or by 
$$X_{m+1}=\frac{1}{Y_{2m+1}+\frac{1}{Y_{2m+2}+X_m}}$$ in the general case.

This characterization by continued fractions is one key ingredient in the proof given by Letac and Weso\l owski~\cite{letw} of the characterization of GIG laws by the Matsumoto-Yor property, which we now recall.

\subsection{The Matsumoto-Yor property}
  Consider two independent, positive random variables $X$ and $Y$ such that $$X \sim {\rm GIG}(-p, a,b), \ \ Y \sim \Gamma(p,a/2)$$
  for some  $p, a, b >0$.
The Matsumoto-Yor property is the fact that the random variables \begin{equation} \label{defuv}
U=\frac{1}{X+Y} \ \ \ \mbox{and} \ \ \ V=\frac{1}{X}-\frac{1}{X+Y}
\end{equation}
 are independent. This property has been discovered in Matsumoto and Yor~\cite{maty1} for the case $a=b$ while studying certain exponential functionals of Brownian motion. Letac and Weso\l owski have noticed afterwards that this property remains true if $a\neq b$, although their paper has finally been published earlier than the Matsumoto-Yor paper from 2001, namely in 2000. Letac and Weso\l owski~\cite{letw} have further proved that this property is in fact a characterization of GIG laws (more exactly, of the product of a GIG and a Gamma law with suitable parameters).

\begin{theor} (Letac and Weso\l owski,~\cite{letw}). \label{Letac-Weso}
Consider two non-Dirac, positive and independent random variables $X$ and $Y$.
Then the random variables $U=(X+Y)^{-1}$ and $V=X^{-1}-(X+Y)^{-1}$ are independent if and only if
there exist $p >0$,   $a>0$ and $b>0$ such that
 $X \sim {\rm GIG}(-p, a,b)$ while  $Y \sim \Gamma(p,a/2)$. \end{theor}\vspace{2mm}

Let us make a few remarks on the theorem:
\begin{itemize}
\item
As shown in the same paper, this result holds true also for matrix variates (characterization of the product of matrix GIG and Wishart variables) under a smoothness assumption not needed in the scalar case.
To prove this, the authors use the extension to the matrix case, established by Bernadac~\cite{ber}, of the
above mentioned  continued-fraction characterization of GIG distributions. 

\item
An interpretation of this property in terms of  Brownian motion has been pointed out by  Matsumoto and Yor~\cite{maty2}. Massam and Weso\l owski~\cite{masw1} have provided a tree-version of the Matsumoto-Yor property, and an interpretation of that tree-version by means of a family of Brownian motions  is given in Weso\l owski and Witkowski~\cite{wesw}.

\item For $p=-1/2$,  the Matsumoto-Yor property can be seen as a consequence of an independence property established by Barndorff-Nielsen and Koudou~\cite{bark} on a tree-network of Inverse Gaussian resistances (see Koudou~\cite{kou}).
 
\item Koudou and Vallois~(\cite{kv1}, \cite{kv}) investigated  a generalization of the Matsumoto-Yor property,
by looking for smooth decreasing functions $f$ from  $(0,\infty)$ onto $(0,\infty)$ with the following property: there exist independent, positive random variables $X$ and $Y$ such that the variables $U=f(X+Y)$ and $V=f(X)-f(X+Y)$ are independent. This led to other independence properties of the Matsumoto-Yor type. In particular, an independence property characterizing Kummer distributions has been established.
\end{itemize}

\subsection{Characterizations by constant regression}
In a subsequent work, Weso\l owski~\cite{wes} relaxed the independence condition between $U$ and $V$ in Theorem \ref{Letac-Weso} by assuming that   $V$ and  $1/V$ have constant regression on $U$. Following Lukacs~\cite{lukacs}, we say that a random variable $Z$ with finite expectation has a constant regression on a random variable $T$ if $\E(Z|T)=\E(Z)$ almost surely. As can be seen in the following theorems, the price to pay for the relaxation at the level of independence is to assume the existence of moments of $X$ and/or $1/X$.

\begin{theor} \label{wes1} (Weso\l owski~\cite{wes}). 
Consider two non-Dirac, positive and independent random variables $X$ and $Y$ such that $\E (1/X)< \infty$.
Let $U$ and $V$ be defined by (\ref{defuv}), and assume that $\E(V|U)=c$ and $\E(1/V|U)=d$ for some real constants $c,d$. \\
Then $cd>1$ and there exists $a>0$ such that  
 $X \sim {\rm GIG}(-p, 2a,2b)$ and $Y \sim \Gamma(p, a)$ where $p=cd/(cd-1)$ and $b=d/(cd-1)$.
\end{theor}\vspace{2mm}

One can prove that  $X$ is GIG distributed,  assuming that  $Y$ is Gamma distributed, and vice-versa :

\begin{theor} \label{wes2} (Seshadri and Weso\l owski~\cite{sesw1}).
Consider two non-Dirac, positive and independent random variables $X$ and $Y$ such that $\E (1/X)$ and $\E (X)$ are finite. Define  $U$ and $V$ as in  (\ref{defuv}).
Suppose $Y \sim \Gamma(p, a/2)$  for some $p>0, a>0$. If the regression of $V$ on $U$ is constant, i.e. $\E(V|U)=c$, then  $X \sim {\rm GIG}(-p, a,b)$, with $b=2p/c$.
\end{theor}\vspace{2mm}

\begin{theor} \label{wes3} (Seshadri and Weso\l owski~\cite{sesw1}).
Consider two non-Dirac, positive and independent random variables $X$ and $Y$ such that $\E (Y)< \infty$.
Assume that $X \sim {\rm GIG}(-p, a,b)$ for some positive $p,a,b$. If $\E(V|U)=2p/b$, then $Y \sim \Gamma(p, a/2)$.
\end{theor}\vspace{2mm}

Theorems \ref{wes1}, \ref{wes2} and \ref{wes3} are established by proving that Laplace transforms of measures linked with the laws of $X$ and $Y$ are probabilistic solutions of second-order differential equations.  Weso\l owski~\cite{wes} has pointed out that the condition $\E(V|U)=c$ can be expressed in terms of linearity of the regression of $Y/X$ on $X+Y$, under the form
$$\E\left(\left.\frac{Y}{X}\, \right|\, X+Y\right)=c (X+Y).$$
This can easily be seen from the definition of $U$ and $V$ in (\ref{defuv}).

The results of Weso\l owski~\cite{wes} and Seshadri and Weso\l owski~\cite{sesw1}    have been slightly extended in 
Chou and Huang~\cite{cho}.
\begin{theor} (Chou and Huang~\cite{cho}).
Consider two non-Dirac, positive and independent random variables $X$ and $Y$ and define $U$ and $V$ as in  (\ref{defuv}).
Assume that, for some fixed integer $r$, $\E (X^{-r-2})$, $\E (X^{-r})$, $\E (Y^r)$ and  $\E (Y^{r+2})$  are finite. \\
Assume that, for some constants $c_r$ and $c_{r+1}$,
$$\E (V^{r+1} |\,U)=c_r \E (V^r |\,U)$$
and
$$\E (V^{r+2} |\,U)=c_{r+1} \E (V^{r+1} |\,U). $$
Then, $c_{r+1}>c_r>0$  and  there exists $a>0$ such that $X \sim {\rm GIG}(-p, 2a,2b)$  and $Y \sim \Gamma(p, a)$ with $p=c_r/ (c_{r+1}-c_r)-r >0$ and
$b=1/ (c_{r+1}-c_r) >0$.
\end{theor}\vspace{2mm}

Lukacs~\cite{lukacs} (Theorem 6.1) characterized the common normal distribution of \mbox{i.i.d.} random variables $X_1,X_2, \ldots, X_n$ by the constancy of regression of a quadratic function of $X_1,X_2, \ldots, X_n$ on the sum  $ \Lambda=\sum_{i=1}^n X_i$. The following theorem characterizing the GIG distribution belongs also to the family of characterizations of distributions by constancy of regression of a statistic $S$ on the sum $\Lambda$ of observations.

\begin{theor} (Pusz~\cite{pus})
Let $X_1,X_2, \ldots, X_n$ be \ind and identically distributed copies of a random variable $X>0$ with
$\E(1/X^2)$,  $\E(1/X)$ and $\E(X)$  finite. Consider $q>0$ and $p\in {\mathbb R}$ such that $p\E(1/X) + q\E(1/X^2)>0$ and define $ \Lambda=\sum_{i=1}^n X_i$ and
$$S=\Lambda \left( p\sum_{i=1}^n \frac{1}{X_i}  +
q\sum_{i=1}^n \frac{1}{X_i^2} \right)-nq\sum_{i=1}^n \frac{1}{X_i}.$$
Then $\E(S|\Lambda)=c$ for some constant $c$ if and only if there exist $\mu \in \R$, $a, b, \delta >0$  such that
$$p=\delta (\mu-1), \, q=\delta b/2, \, c=\delta n(n\mu-1)$$
and $X \sim {\rm GIG} (p, a, b)$.
\end{theor}\vspace{2mm}
The proof is based on the fact that the GIG distribution is characterized by its moments, on a lemma by Kagan {\it et al.}~\cite{kag} giving a necessary and sufficient condition for the constancy of $\E(S|\Lambda)$, and on a  careful  manipulation of a differential equation satisfied by the function 
$f(t)= \E(1/X^2)\exp(itX)$.

\subsection{Entropy characterization}

 Characterizing
a distribution by the maximum entropy principle  dates back to Shannon~\cite{sha} who showed 
that Gaussian random variables maximize entropy among all real-valued random variables with given mean and variance. 
Since then,  many examples of such  characterizations have appeared in the literature. For instance, Kagan {\it et al.}~\cite{kag} { characterized several well-known distributions such as the exponential, Gamma or Beta  distributions in terms of maximum entropy given  various constraints, e.g.:
 {\it (i)} the exponential  distribution with parameter $\lambda$  maximizes the entropy among all the distributions of $X$ supported on the positive real line under the constraint $\E (X)= 1/\lambda$;  {\it (ii)} the Gamma distribution $\Gamma (p, a)$ maximizes the entropy under the constraints $\E (X)= p/a$ and $\E (\log X)= \Psi (p)- \log a$, where $\Psi$ is the Digamma function; {\it (iii)} the Beta distribution with parameters $a, b>0$ maximizes the entropy among all the distributions of $X$ supported on $[0,1]$ under the constraints $\E (\log X)= \Psi (a)- \Psi (a+b)$ and $\E (\log (1-X))= \Psi (b)- \Psi (a+b)$.}   
 This type of \crt exists as well for GIG distributions, { as established in the following theorem due to Kawamura and K${\bar {\rm o}}$sei~\cite{kaw}, for which we recall that the Bessel function of the third kind is defined by 
$$K_p (z)=2^{-p -1}z^{p} \int_0^\infty x^{-p -1} e^{-x-\frac{z^2}{4x}} dx, \ \ {\rm Re} (z)>0.$$}
\begin{theor} \label{Kawamura}
The distribution of $X$ with density $f$ on $(0,\infty)$ which maximizes the entropy
$$ H(X)=-\int_0^\infty f(x) \log f(x)\, dx$$ 
under the constraints 
$$\E  \left(  \log \frac{X}{\sqrt{b/a}}      \right)=\partial _p \log K_p(\sqrt{ab})\,\mbox{and}\,
\ \ \E \left(   \frac{X}{\sqrt{b/a}} +    \frac{\sqrt{b/a}}{X}  \right)=
\frac{K_{p+1}(\sqrt{ab})+   K_{p-1}(\sqrt{ab})     }{K_p(\sqrt{ab})} $$
is the distribution $GIG(p,a,b)$.
\end{theor}\vspace{2mm}

One application of entropy characterizations are goodness-of-fit tests. An instance is Vasicek~\cite{vas} who designed an entropy test of normality, using Shannon's characterization. In our GIG setting, Mudholkar and Tian~\cite{mudt} established a goodness-of-fit test for the special case of an
Inverse Gaussian (IG) distribution based on an entropy characterization that we recall below.
Since the IG is a major particular case of GIG distributions ---the latter were actually introduced in order to meet the needs of adding a third parameter to the IG distribution---we devote a brief subsection to characterizations of the IG distribution.

\subsection{Some characterizations of the IG distribution}
Of particular interest for the present survey is the Inverse Gaussian  distribution with density
$$\frac{\sqrt{b}}{\sqrt{2\pi}}e^{\sqrt{ab}}x^{-\frac{3}{2}} 
e^{-(ax+b/x)/2}, \ \ \ x >0.$$
The IG distribution is useful  in many fields  for data modeling and analysis when
the observations are right-skewed. Numerous examples of applications can be found in 
Chhikara and Folks~\cite{chh} and Seshadri~\cite{Ses99}.
The IG distribution possesses many similarities in terms of statistical properties with the Gaussian distribution, as pointed out by \mbox{e.g.}  Mudholkar and Tian~\cite{mudt}. We give here a few results characterizing  the IG distribution. 
\subsubsection{Entropy characterizations}
Let us first observe that,
as an immediate  corollary of Theorem \ref{Kawamura}, we obtain the following (new) entropy characterization of the IG distribution.
\begin{theor}
Denote by $c_{a,b}$ the value of $\partial _p \log K_p(\sqrt{ab})$ at $p=-1/2$.
The distribution of $X$ with density $f$ on $(0,\infty)$ having maximum entropy
under the constraints 
$$\E  \left(  \log \frac{X}{\sqrt{b/a}}      \right)=c_{a,b}\ \ \,\mbox{and}\,\,
\ \ \E \left(   \frac{X}{\sqrt{b/a}} +    \frac{\sqrt{b/a}}{X}  \right)=
2+\frac{1}{\sqrt{ab}}$$
is the distribution ${\rm GIG}(-1/2,a,b)={\rm IG}(a,b)$.
\end{theor}\vspace{2mm}

 Mudholkar and Tian~\cite{mudt} have obtained another entropy characterization of the IG distribution.

\begin{theor} ( Mudholkar and Tian~\cite{mudt}).
A random variable $X$ has the IG$(\sqrt{b/a}, b)$ distribution if and only if $1/\sqrt{X}$ attains maximum entropy among all
absolutely continuous random variables $Y\geq 0$ such that
$$\E (Y^{-2})=\sqrt{b/a}\ \ \,\,\mbox{and} \, \, \ \ \E (Y^{2})=\sqrt{a/b}+1/b.$$
\end{theor}\vspace{2mm}
{Inspired by the Vasicek~\cite{vas} methodology for the entropy test of normality,  Mudholkar and Tian~\cite{mudt}, as stated earlier, used the above theorem to construct a goodness-of-fit test for the IG distribution. The test statistic is based on a non-parametric estimation of the entropy of $1/\sqrt{X}$ and rejection occurs when this estimation is not close enough to the theoretical value  under the null hypothesis of $X$ being IG-distributed. The rejection criterion is assessed at given significance level by critical values obtained via Monte Carlo simulations.} 

\subsubsection{A martingale characterization}
Seshadri and Weso\l owski~\cite{sesw4} established the following martingale characterization of the IG distribution.
\begin{theor} \label{mart}(Seshadri and Weso\l owski~\cite{sesw4}).
Let $(X_n)_{n\geq 1}$ be a sequence of positive, non degenerate  \mbox{i.i.d.} random variables. For $n\geq 1$, define $S_n= X_1+X_2+\cdots+ X_n$ and consider the $\sigma$-algebra $\mathcal{F}_n= \sigma(S_n, S_{n+1}, \ldots)$. Let $b>0$. The sequence   
$$\left(  \frac{n}{S_n}-\frac{1}{2bn}, \,  \mathcal{F}_n \right)_{n\geq 1}$$
is a backward martingale if and only if $X_1 \sim {\rm IG} (a, b)$ for some $a>0$.
\end{theor}\vspace{2mm}

{\bf Remark:} In fact, the result established in Seshadri and Weso\l owski~\cite{sesw4} is more general than stated in Theorem \ref{mart}. They proved that $(  \frac{\alpha_n}{S_n}-\beta_n, \,  \mathcal{F}_n )_{n\geq 1}$ is a backward martingale for some sequences $(\alpha_n)_n$ and $(\beta_n)_n$ if and only if
$X_1$ follows one of four distributions, among which the IG distribution. The proof uses the Letac and Mora~\cite{letm} \crt of natural exponential families with cubic variance function.

We conclude this section by writing down the well-known Khatri characterization of the IG distribution.

\begin{theor}
(Khatri~\cite{kha}). Let $X_1,X_2, \ldots, X_n$ be \mbox{i.i.d.} copies of a random variable $X$ such that
$\E(1/X)$,  $\E(X)$, $\E(X^2)$ and $\E(1/\sum_{i=1}^n X_i)$ exist and do not vanish.\\
Let $$\bar{X} = \frac{1}{n}\sum_{i=1}^n X_i \ \ and \ \ \bar{X}_{-1} = \frac{1}{n}\sum_{i=1}^n \frac{1}{X_i}.$$
The random variables $\bar{X}$  and $ \bar{X}_{-1}-1 / \bar{X}$
are independent if and only if there exist $a,b>0$ such that $ X\sim {\rm IG}(a,b)$.
\end{theor}\vspace{2mm}

If one considers the parameters $\mu=\sqrt{b/a}$ and $\lambda=b$, then 
$\bar{X}$ and $1/(\bar{X}_{-1}-1/\bar{X})$ are respectively the  maximum likelihood estimators of $\mu$ and $\lambda$. The independence between these estimators characterizes the IG distribution, as the independence between the empirical mean and standard deviation as maximum likelihood estimators of location and scale characterizes
the normal distribution.  This is one of the similarities between the IG and normal distributions, as observed  by Mudholkar and Tian~\cite{mudt}.

\section{New characterizations of GIG laws}
\subsection{An MLE characterization}

A famous characterization theorem in statistics  due to Carl Friedrich Gauss~\cite{Gauss}  says the following:  the sample mean $\bar{X}$ is for any \mbox{i.i.d.} sample $X_1,\ldots,X_n$ of any sample size~$n$ the maximum likelihood estimator (MLE) of the location parameter in a location family $\{f(x-\mu), \mu\in\R\}$ if and only if the samples are drawn from a Gaussian population (with variance not fixed). This very first \emph{MLE characterization theorem} has important implications, as it clearly indicates that, the further one is away from the Gaussian situation, the less efficient becomes the sample mean as estimator. Several other characterization theorems have emerged in the past, linking particular forms of MLEs to a specific distribution such as the one-parameter exponential family (Poincar\'e~\cite{Poincare}), the (negative) exponential distribution (Teicher~\cite{Tei}), the Laplace distribution (Kagan \emph{et al.}~\cite{kag}), the Gamma distribution (Marshall and Olkin~\cite{MO}) or the Harmonic Law (H\"urlimann~\cite{Hur}). For a recent overview and a unified theory on these results, we refer to Duerinckx \emph{et al.}~\cite{due}. 

In this section, our aim is to apply the general result of Duerinckx \emph{et al.}~\cite{due} to the GIG distribution in order to construct an MLE characterization theorem for GIG laws. To this end, we shall rather have recourse to the re-formulation~(\ref{GIG2}) of the density, as in that parameterization the family $\{ f_{p,\theta, \eta}(x),\,\eta >0\}$ is a scale family. Let us first observe, by a simple calculation, that the MLE    of  $\eta$ for  fixed $p$ and $\theta$ is
\begin{equation*}
\hat{\eta}= \frac{\sqrt{p^2+\theta^2 \bar{X}\bar{X}_{-1}} -p}{\theta  \bar{X}_{-1} } , 
\end{equation*}
where we recall that $\bar{X} = \frac{1}{n}\sum_{i=1}^n X_i$  and $\bar{X}_{-1} = \frac{1}{n}\sum_{i=1}^n \frac{1}{X_i}$
if $X_1,\ldots X_n$ are \mbox{i.i.d.} with common distribution ${\rm GIG} (p, a, b)$. This enables us to formulate the following.

\begin{theor}\label{MLE}
Consider an integer $n\geq 3$.
For $p\in \R$ and $\theta >0$, $\frac{\sqrt{p^2+\theta^2 \bar{X}\bar{X}_{-1}} -p}{\theta  \bar{X}_{-1} } $ is the MLE  of the scale parameter $\eta$ in a scale family $\left\{ \frac{1}{\eta} f(x/ \eta), \eta >0\right\}$ on $(0,\infty)$  for all samples $X_1,\ldots,X_n$ of fixed sample size $n$ if and only if there exists $d>0$ such that the family is  $\{ f_{pd, \theta d, \eta}(x), \, \eta >0\}$, i.e. the family of densities of the distributions ${\rm GIG} (pd, \theta d/\eta, \theta \eta d)$, $\eta >0$.
\end{theor}\vspace{2mm}

\noindent{\sc Proof of Theorem~\ref{MLE}.} We want to apply Theorem~5.1 of Duerinckx \emph{et al.}~\cite{due}. To this end, we first need to calculate the scale score function of the standardized GIG density $f_{p,\theta}(x):=\frac{1}{ 2K_p (\theta)} x^{p -1}e^{-\frac{1}{2}\theta( x+  1/x)}$ defined on $(0, \infty)$. This score, defined as $\psi_{f_{p,\theta}}(x)=1+x\frac{f'_{p,\theta}(x)}{f_{p,\theta}(x)}$, corresponds to
$$\psi_{f_{p,\theta}}(x)=p-\frac{\theta}{2}x+ \frac{\theta}{2x}.$$
One easily sees that $\psi_{f_{p,\theta}}$ is invertible over $(0,\infty)$ with image  $\R$. Thus Theorem~5.1 of Duerinckx \emph{et al.}~\cite{due} applies and the so-called MNSS (Minimal Necessary Sample Size, the smallest sample size for which an MLE characterization holds) equals 3. This theorem then yields  that $\frac{\sqrt{p^2+\theta^2 \bar{X}\bar{X}_{-1}} -p}{\theta  \bar{X}_{-1} } $ is the MLE  of the scale parameter $\eta$ in a scale family $\left\{ \frac{1}{\eta} f(x/ \eta), \eta >0\right\}$ on $(0,\infty)$  for all samples $X_1,\ldots,X_n$ of fixed sample size $n\geq 3$ if and only if the densities are proportional to $x^{d-1}(f_{p,\theta}(x))^d$ for $d>0$, which is nothing else but $f_{pd,\theta d}(x)$, the (standardized) GIG density with parameters $pd$ and $\theta d$. This yields the announced MLE characterization for GIG laws. \cqfd \vspace{2mm}

Clearly, this MLE characterization theorem for GIG laws contains, \emph{inter alia}, MLE characterizations for the Inverse Gaussian, the Reciprocal Inverse Gaussian and the hyperbolic or Harmonic Law, hence the characterization of H\"urlimann~\cite{Hur} is a special case of our Theorem~\ref{MLE}.

\subsection{A Stein characterization}
The celebrated Stein's method of normal approximation, introduced in Stein~\cite{ste}, has over the years been adapted to several other probability distributions, including the Poisson (Chen~\cite{che}), the exponential (Chatterjee {\it et al.}~\cite{cha}), the Gamma (Luk~\cite{luk}), the multinomial (Loh~\cite{loh}), the geometric (Pek\"oz~\cite{pek}), the negative binomial (Brown and Phillips~\cite{bro}) or the Beta (Goldstein and Reinert~\cite{gol}). 
A first step in this method consists in finding a suitable Stein operator, whose properties determine the quality of the approximation; see Ross~\cite{ros} for a recent overview on the intricacies of Stein's method. This operator satisfies a \emph{Stein characterization theorem} which clearly links it to the targeted distribution. 

%

To the best of the authors' knowledge, so far no such Stein characterization has been proposed in the literature for the GIG distribution, whence the reason of the following result. 

\begin{theor}\label{Steintheo}
A positive random variable $X$ follows the GIG$(p,a,b)$ distribution if and only if for any differentiable function $h$
satisfying 
\begin{equation*}
\lim_{x\to \infty} f_{p,a,b}(x) h(x) =\lim_{x\to 0} f_{p,a,b}(x) h(x)=0
\end{equation*}
we have
\begin{equation}\label{Steinchar}
\E \left[h^\prime (X) + \left( \frac{p-1}{X} 
+\frac{b}{2X^2} -\frac{a}{2} \right) h(X)  \right]=0.
\end{equation}
The functional $h\mapsto \mathcal{T}_{f_{p,a,b}}(h)(x):=h^\prime (x) + \left( \frac{p-1}{x}+\frac{b}{2x^2} -\frac{a}{2} \right) h(x)$ is the GIG Stein operator. 
\end{theor}\vspace{2mm}

\noindent{\sc Proof of Theorem~\ref{Steintheo}.} The sufficient condition is readily checked by noting that $h^\prime (x) + \left( \frac{p-1}{x} 
+\frac{b}{2x^2} -\frac{a}{2} \right) h(x)=\frac{(h(x)f_{p,a,b}(x))'}{f_{p,a,b}(x)}$, since then
$$
\E \left[h' (X) + \left( \frac{p-1}{X} 
+\frac{b}{2X^2} -\frac{a}{2} \right) h(X)  \right]=\int_0^\infty (h(x)f_{p,a,b}(x))' dx=0
$$
by the conditions on $h$. To see the necessity, write $f$ for the density of $X$  and define for $z>0$ the function  
\begin{equation*}
  l_z(u):= {\mathbb{I}}_{(0, z]}(u) -F_{p,a,b}(z), \quad u>0,
\end{equation*}
 with $F_{p,a,b}$ the GIG cumulative distribution function and $\mathbb{I}$ an indicator function. Then the function 
$$
    h_z(x) : = \frac{1}{f_{p,a,b}(x)} \int_{0}^x l_z(u) f_{p,a,b}(u) du \left( =
      -\frac{1}{f_{p,a,b}(x)} \int_{x}^\infty l_z(u) f_{p,a,b}(u) du \right) 
$$
is differentiable and satisfies $ \lim_{x\to \infty} f_{p,a,b}(x) h_z(x) =\lim_{x\to 0} f_{p,a,b}(x) h_z(x)=0$ for all $z$ (since $\int_0^\infty l_z(u)f_{p,a,b}(u)=0$), hence is a candidate function for the functions $h$ that verify~\eqref{Steinchar}. Using hence this $h_z$ leads to 
\begin{eqnarray*}
0&=&\E \left[h_z' (X) + \left( \frac{p-1}{X} +\frac{b}{2X^2} -\frac{a}{2} \right) h_z(X)  \right]\\
&=&\E\left[\frac{(h_z(X)f_{p,a,b}(X))'}{f_{p,a,b}(X)}\right]\\
&=&\int_0^\infty l_z(x) f(x)dx\\
&=& F(z)-F_{p,a,b}(z)
\end{eqnarray*}
for every $z>0$, from which we can conclude that $F=F_{p,a,b}$. \cqfd\vspace{2mm}

This result is a particular instance of the \emph{density approach} to Stein characterizations initiated in Stein \emph{et al.}~\cite{Steinetal} and further developed in Ley and Swan~\cite{leys}. We attract the reader's attention to the fact that we could replace the functions $h(x)$ with $h(x)x^2$, turning the GIG Stein operator $\mathcal{T}_{f_{p,a,b}}(h)$ into the perhaps more tractable form
$$
\mathcal{T}_{f_{p,a,b}}(h)(x)=x^2 h^\prime (x) + \left(-x^2\frac{a}{2} + (p+1)x +\frac{b}{2} \right) h(x).
$$
This new Stein characterization and the associated Stein operator(s) could be exploited to derive rates of convergence in some asymptotics related to GIG distributions, e.g. the rate of convergence of the law of the continued fraction in (\ref{contf}) to the GIG distribution.    This kind of application will be investigated in future work.

\vspace{1cm}

\noindent ACKNOWLEDGEMENTS:

\noindent This research has been supported by Wallonie-Bruxelles International, the Fonds de la Recherche Scientifique, the Minist\`ere Fran\c cais des Affaires \'etrang\`eres et europ\'eennes and the Minist\`ere de l'Enseignement sup\'erieur et de la Recherche as part of the Hubert Curien 2013 Grant. Christophe Ley, who is also a member of ECARES,   thanks the Fonds National de la Recherche Scientifique, Communaut\'e Fran\c caise de Belgique, for support via a Mandat de Charg\'e de Recherche FNRS.

\end{document}